\documentclass[12pt,a4paper]{amsart}
\usepackage{url,amscd,amssymb}

\newcommand{\invlim}{\lim_{\leftarrow}\limits}
\newcommand{\Q}{\ensuremath{\mathbb{Q}}}
\newcommand{\Z}{\ensuremath{\mathbb{Z}}}
\newcommand{\N}{\ensuremath{\mathbb{N}}}
\newcommand{\C}{\ensuremath{\mathbb{C}}}
\newcommand{\Zp}{\ensuremath{\Z_{p}}}

\newcommand{\Zpl}{\ensuremath{\Z_{(p)}}}
\newcommand{\Zplu}{\ensuremath{\Zpl^{\times}}}
\newcommand{\Zt}{\ensuremath{\Z_{2}}}
\newcommand{\Ztu}{\ensuremath{\Zt^{\times}}}
\newcommand{\Ztl}{\ensuremath{\Z_{(2)}}}

\newcommand{\padic}{$p$\nobreakdash-adic}
\newcommand{\padically}{\padic ally}
\newcommand{\plocal}{$p$\nobreakdash-local}
\newcommand{\twolocal}{$2$\nobreakdash-local}
\newcommand{\Ktheory}{$K$\nobreakdash-theory}
\newcommand{\KOtheory}{$KO$\nobreakdash-theory}
\newcommand{\floor}[1]{\left\lfloor#1\right\rfloor}
\newcommand{\ceiling}[1]{\left\lceil#1\right\rceil}

\newcommand{\gp}[2]{\genfrac{[}{]}{0pt}{}{#1}{#2}}

\newcommand{\KzK}{\ensuremath{K_{0}(K)}}

\newcommand{\Kzk}{\ensuremath{K_{0}(k)}}
\newcommand{\KuzK}{\ensuremath{K^{0}(K)}}
\newcommand{\KusK}{\ensuremath{K^{*}(K)}}
\newcommand{\Kuzk}{\ensuremath{K^{0}(k)}}
\newcommand{\kuzk}{\ensuremath{k^{0}(k)}}
\newcommand{\kusk}{\ensuremath{k^{*}(k)}}
\newcommand{\guzg}{\ensuremath{g^{0}(g)}}
\newcommand{\Guzg}{\ensuremath{G^{0}(g)}}
\newcommand{\GuzG}{\ensuremath{G^{0}(G)}}
\newcommand{\Gzg}{\ensuremath{G_{0}(g)}}
\newcommand{\tensor}{\otimes}
\newcommand{\set}[1]{\ensuremath{\{\,#1\,\}}}

\newcommand{\CP}{\ensuremath{\C P^{\infty}}}

\newcommand{\At}{\ensuremath{\widetilde{A}}}

\let\bigset\setbig

\let\le\leqslant

\let\ge\geqslant
\let\setminus\smallsetminus
\newcommand{\bone}{\ensuremath{{\mathbf 1}}}
\newcommand{\ba}{\ensuremath{\mathbf{a}}}
\newcommand{\bb}{\ensuremath{\mathbf{b}}}
\newcommand{\bc}{\ensuremath{\mathbf{c}}}
\newcommand{\bq}{\ensuremath{\mathbf{q}}}
\newcommand{\bqh}{\ensuremath{\mathbf{\qh}}}
\newcommand{\qb}{\ensuremath{\bar{q}}}
\newcommand{\bqb}{\ensuremath{\mathbf{\qb}}}
\newcommand{\fh}{\hat{f}}
\newcommand{\phih}{\hat{\phi}}
\newcommand{\thetah}{\hat{\theta}}
\newcommand{\thetab}{\bar{\theta}}
\newcommand{\qh}{\hat{q}}

\let\phi\varphi
\let\epsilon\varepsilon

\newtheorem{thm}[equation]{Theorem}
\newtheorem{prop}[equation]{Proposition}
\newtheorem{cor}[equation]{Corollary}
\newtheorem{lemma}[equation]{Lemma}

\theoremstyle{definition}
\newtheorem{defn}[equation]{Definition}
\newtheorem{remark}[equation]{Remark}

\numberwithin{equation}{section}

\begin{document}

\title{Algebras of operations in $K$-theory}

\author[F. Clarke]{Francis Clarke}
\address{Department of Mathematics, University of Wales Swansea,
Swansea SA2 8PP, Wales}
\email{F.Clarke@Swansea.ac.uk}

\author[M. D. Crossley]{Martin Crossley}
\address{Department of Mathematics, University of Wales Swansea,
Swansea SA2 8PP, Wales}
\email{M.D.Crossley@Swansea.ac.uk}

\author[S. Whitehouse]{Sarah Whitehouse}
\address{Department of Pure Mathematics,
University of Sheffield,
Sheffield S3 7RH, England}
\email{S.Whitehouse@sheffield.ac.uk}

\begin{abstract}
    We describe explicitly the algebras of degree zero operations in
    connective and periodic $p$-local complex $K$-theory.  Operations
    are written uniquely in terms of certain infinite linear
    combinations of Adams operations, and we give formulas for the
    product and coproduct structure maps.  It is shown that these
    rings of operations are not Noetherian.  Versions of the results
    are provided for the Adams summand and for real $K$-theory.
\end{abstract}

\keywords{\Ktheory\ operations, Gaussian polynomials}

\subjclass[2000]{%
Primary:   55S25; 
Secondary: 19L64, 
           11B65. 
}

\date{$21^{\text{st}}$ January 2004}

\maketitle

\tableofcontents

\section{Introduction}
\label{intro}

The complex \Ktheory\ of a space or spectrum may be usefully endowed
with operations.  The most well-known of these are the Adams
operations $\Psi^{j}$ arising out of the geometry of vector bundles. 
As originally constructed by Adams, these are {\it unstable}
operations.  A stable operation is given by a sequence of maps,
commuting with the Bott periodicity isomorphism.  For an Adams
operation $\Psi^{j}$ to be {\it stable} requires $j$ to be a unit in
the coefficient ring one is working over.  Integrally, we only have
$\Psi^{1}$ and $\Psi^{-1}$, corresponding to the identity and complex
conjugation.  Following work of Adams, Harris and Switzer~\cite{ahs}
in which the structure of the dual object, the algebra of
cooperations, was determined, Adams and Clarke~\cite{ac} showed that
there are uncountably many integral stable operations.  It is
remarkable that, after more than forty years of topological \Ktheory,
no-one knows explicitly any integral stable operations, apart from
linear combinations of the identity and complex conjugation.

In this paper we give a new description of \Ktheory\ operations in the
\plocal\ setting.  This is considerably closer to the integral
situation and more subtle than the rather better understood case of
\padic\ coefficients; see~\cite{Clarke1987,mst}.  In \plocal\
connective \Ktheory, the degree zero operations $k_{(p)}^{0}(k_{(p)})$
form a bicommutative bialgebra, which we will denote by~$\kuzk_{(p)}$. 
Note, however that this is not the $p$\nobreakdash-localisation of the
bialgebra~\kuzk\ of integral operations, but is isomorphic to the
\emph{completed} tensor product $\kuzk\widehat{\otimes}\Zpl$.  In the
periodic theory, the corresponding object
$\KuzK_{(p)}=K_{(p)}^{0}(K_{(p)})=\KuzK\widehat{\otimes}\Zpl$ is a
bicommutative Hopf algebra, as it also possesses an antipode.  We
provide explicit descriptions for both of these algebras of
operations, together with formulas for all their structure maps. 
These results build on our recent paper~\cite{ccw}, in which we gave
new additive bases for the ring of cooperations in \plocal\ \Ktheory. 
Such an understanding of the degree zero operations is sufficient to
determine the whole of $\KusK_{(p)}$ and the torsion-free part of
$\kusk_{(p)}$; see~\cite{Lellman1984}.

Since we have a stable Adams operation $\Psi^{\alpha}$ for each
$\alpha\in \Zplu$, and since these multiply according to the formula
$\Psi^{\alpha}\Psi^{\beta}=\Psi^{\alpha\beta}$, the group ring
$\Zpl[\Zplu]$ is a subring of the ring of operations.  Our results
express operations in terms of certain infinite sums involving Adams
operations and thus describe the ring of operations as a completion of
the group ring.  This idea is certainly not new;
Johnson~\cite{Johnson1986} also has basis elements for
\plocal\ operations of this form.  However, our description has
considerable advantages in the form of explicit formulas, which in the
connective case are particularly nice.  We note that Madsen, Snaith
and Tornehave~\cite{mst} have also considered operations defined as
infinite sums of Adams operations, but for them the \padic\ context is
essential.

Our results obtained will be used in a later paper to give a
simplification of Bousfield's work in~\cite{Bousfield85} describing
the $K_{(p)}$\nobreakdash-local category.

We now outline the structure of the paper.

Sections~\ref{hopfalgcon} and~\ref{ringstructure} are concerned with the 
case of operations in connective \Ktheory.  In Section~\ref{hopfalgcon} 
we describe the bialgebra~$\kuzk_{(p)}$ in an explicit form which 
enables us to also describe the structure maps.  We also give formulas for 
the action on the coefficient ring and on the Hopf bundle over~\CP.  
In Section~\ref{ringstructure} we show that, as a ring, $\kuzk_{(p)}$ 
is not Noetherian, and we characterise its units.  We also indicate 
how it can be considered as a completion of a polynomial ring.

In Section~\ref{adamscon} we consider the idempotents in connective
\Ktheory\ which were introduced by Adams, and we show how the results
of Sections~\ref{hopfalgcon} and~\ref{ringstructure} extend to the
algebra of operations on the Adams summand.  We prove in
Section~\ref{pcomp} that the ring of operations on the \padic\ Adams
summand is a power series ring.

In Section~\ref{hopfalgper} we show how the results of
Sections~\ref{hopfalgcon}--\ref{adamscon} generalise to periodic 
\Ktheory, and in Section~\ref{convper} we discuss the relation 
between the connective and periodic cases.

In Sections~\ref{2-local} and~\ref{KO} we work over the prime~$2$.  We 
outline how our results from the preceding sections need to be adapted, 
and we consider operations in \KOtheory.

Finally, in an appendix we give some general relations among
polynomials which underpin a number of the formulas given in the
preceding sections.

Unless otherwise stated, $p$ is assumed to be an odd prime.  Having
chosen $p$, we fix $q$ to be an integer which is primitive
modulo~$p^{2}$, and hence primitive modulo any power of~$p$.

\medskip

The third author acknowledges the support of a Scheme~4 grant from 
the London Mathematical Society.

\section{Degree zero operations in connective $K$-theory}
\label{hopfalgcon}

For each non-negative integer $n$, we define $\theta_{n}(X)\in\Z[X]$
by
$$
\theta_{n}(X)
=
\prod_{i=0}^{n-1} (X-q^{i}),
$$
where $q$, as stated in the Introduction, is primitive modulo~$p^{2}$. 
The notation derives from~\cite{Ihrig1981}.  Generalisations of these
polynomials are considered later in this paper.

The Gaussian polynomials in the variable~$q$ (or
$q$\nobreakdash-binomial coefficients) may be defined as
$$
\gp{n}{j}
=
\frac{\theta_{j}(q^{n})}{\theta_{j}(q^{j})}.
$$

\begin{defn}
    Define elements $\phi_{n}\in\kuzk_{(p)}$, for $n\ge 0$, by
    $$
    \phi_{n}
    =
    \theta_{n}(\Psi^{q}),
    $$
    where $\Psi^{q}$ is the Adams operation.
\end{defn}

Thus, for example, $\phi_{0}=1$, $\phi_{1}=\Psi^{q}-1$ and
$\phi_{2}=(\Psi^{q}-1)(\Psi^{q}-q)$.

\begin{thm}\label{phibasis}
    The elements of $\kuzk_{(p)}$ can be expressed uniquely as
    infinite sums
    $$
    \sum_{n\ge 0}a_{n}\phi_{n},
    $$
    where $a_{n}\in\Z_{(p)}$.
\end{thm}

\begin{proof}
   The bialgebra $\kuzk_{(p)}=\Kuzk_{(p)}$ is the $\Z_{(p)}$-dual of
   the bialgebra $\Kzk_{(p)}$; see~\cite{Clarke1987,Johnson1986}.  In
   Proposition~3 of~\cite{ccw} we gave a basis for $\Kzk_{(p)}$,
   consisting of the polynomials
   $f_{n}(w)=\theta_{n}(w)/\theta_{n}(q^{n})$, for $n\ge 0$.  The
   theorem follows from the fact (which is implicit
   in~\cite{Ihrig1981}) that the~$\phi_{n}$ are dual to this basis. 
   To see this, recall that a coalgebra admits an action of its dual
   which, in the case of $\Kzk_{(p)}$, is determined by $\Psi^{r}\cdot
   f(w)=f(rw)$.  By a simple induction on~$m$, this implies that
    $$
    \phi_{m}\cdot f_{n}(w)
    =
    q^{m(m-n)}w^{m}f_{n-m}(w),
    $$
    where $f_{i}(w)$ is understood as $0$ if $i<0$, so that
    $\phi_{m}\cdot f_{n}(w)=0$ if $n<m$.  The Kronecker pairing can
    be recovered from this action by evaluating at $w=1$, hence
    \begin{equation*}
        \bigl\langle\phi_{m},f_{n}(w)\bigr\rangle
        =
        q^{m(m-n)}f_{n-m}(1)
        =
        \begin{cases}
            1,& \text{if $m=n$,}\\ 
            0,& \text{otherwise.} 
        \end{cases}
        \qedhere
    \end{equation*}
\end{proof}

\begin{remark}\label{coeff-action}
    Operations in \Ktheory\ are determined by their action on
    coefficients~\cite{Johnson1984}.  It is therefore instructive to see
    how an infinite sum $\sum_{n\ge 0}a_{n}\phi_{n}$ acts.  Since
    $\Psi^{q}$ acts on $\pi_{2i}(k_{(p)})$ as multiplication
    by~$q^{i}$, we see that $\sum_{n\ge 0}a_{n}\phi_{n}$ acts on the
    coefficient group~$\pi_{2i}(k_{(p)})$ as multiplication by
    $$
    \sum_{n=0}^{i}a_{n}\theta_{n}(q^{i}).
    $$  
    The sum is finite since $\theta_{n}(q^{i})=0$ for $n>i$.
    
    In particular, the augmentation $\epsilon:\kuzk_{(p)}\to\Zpl$ 
    given by the action on~$\pi_{0}(k_{(p)})$, satisfies 
    $\epsilon\left(\sum_{n\ge 0}a_{n}\phi_{n}\right)=a_{0}$.
\end{remark}

It is easy to see by induction that
$\theta_{n}(X)=\sum_{j=0}^{n}(-1)^{n-j}q^{\binom{n-j}{2}}\gp{n}{j}X^{j}$;
see~\cite[(3.3.6)]{Andrews1976} and also \cite[Proposition~8]{ccw}. 
Hence we can express each $\phi_{n}$ explicitly as a linear
combination of Adams operations.

\begin{prop}\label{phi_n-exp}
    For all $n\ge 0$,
    \begin{equation*}
        \phi_{n}
        =
        \sum_{j=0}^{n}(-1)^{n-j}q^{\binom{n-j}{2}}\gp{n}{j}\Psi^{q^{j}}.
    \end{equation*}
    \qed
\end{prop}

Conversely, our proof of Theorem~\ref{phibasis} shows how to express
all the stable Adams operations in terms of the~$\phi_{n}$.

\begin{prop}\label{Psi_j-exp}
    If $j\in\Zplu$,
    $$
    \Psi^{j}
    =
    \sum_{n\ge 0}\frac{\theta_{n}(j)}{\theta_{n}(q^{n})}\phi_{n}.
    $$
    In particular, for $i\in\Z$,
    $$
    \Psi^{q^{i}}
    =
    \sum_{n\ge 0}\gp{i}{n}\phi_{n}.
    $$
    Note that this is a finite sum for $i\ge 0$.
    \qed
\end{prop}

Additive operations in \Ktheory\ are determined by their action on the
Hopf bundle over~\CP; see~\cite{Handbook}.  Writing
$$
k^{0}(\CP)_{(p)}
=
K^{0}(\CP)_{(p)}
=
\Zpl[[t]], 
$$
where $1+t$ is the Hopf bundle, we have the following formula for the
action of~$\kuzk_{(p)}$ on the Hopf bundle.

\begin{prop}
    For all $n\ge 0$,
    $$
    \phi_{n}(1+t)
    =
    \sum_{i\ge n}
    \left(
    \sum_{j=0}^{n}(-1)^{n-j}q^{\binom{n-j}{2}}\binom{q^{j}}{i}\gp{n}{j}
    \right)
    t^{i}.
    $$
\end{prop}

\begin{proof}
    Since $\Psi^{r}$ acts on line bundles by raising them to the
    $r$\nobreakdash-th power, Proposition~\ref{phi_n-exp} shows that
    $$
    \phi_{n}(1+t)
    =
    \sum_{j=0}^{n}(-1)^{n-j}q^{\binom{n-j}{2}}\gp{n}{j}(1+t)^{q^{j}}.
    $$
    The formula now follows by using the binomial expansion and
    reversing the order of summation.  That the coefficient of~$t^{i}$
    in $\phi_{n}(1+t)$ is zero for $i<n$ can be proved by a simple
    induction, using the identity
    $\phi_{n+1}=(\Psi^{q}-q^{n})\phi_{n}$.
\end{proof}

The product structure on $\kuzk_{(p)}$ is determined by the
following formula.

\begin{prop}\label{phi-prod}
    $$
    \phi_{r}\phi_{s}
    = 
    \sum_{i=0}^{\min(r,s)}
    \frac{\theta_{i}(q^{r})\theta_{i}(q^{s})}{\theta_{i}(q^{i})}\,
    \phi_{r+s-i}
    = 
    \sum_{j=\max(r,s)}^{r+s}
    \frac{\theta_{r+s-j}(q^{r})\theta_{r+s-j}(q^{s})}
    {\theta_{r+s-j}(q^{r+s-j})}\,
    \phi_{j}.
    $$
\end{prop}

\begin{proof}
    This result is essentially a fact about the polynomials
    $\theta_{n}(X)$.  Since the algebra of these polynomials underlies
    many of our results, we have gathered together the relevant facts,
    in appropriate generality, in an appendix (Section~\ref{polyid}). 
    In particular, we need to use here Proposition~\ref{prodexp} with
    $m=0$, $X=\Psi^{q}$ and $\bc=\bq=(1,q,q^{2},\dotsc)$.

    To show that the coefficient $A_{r,j-s}(\bq,\bq[s])$ given by that
    proposition is equal to
    $\theta_{r+s-j}(q^{r})\theta_{r+s-j}(q^{s})/\theta_{r+s-j}(q^{r+s-j})$,
    it is only necessary to verify that this expression satisfies the
    recurrence~\eqref{prodrecur}, i.e., taking $i=r+s-j$ that
    $$
    \frac{\theta_{i+1}(q^{r+1})\theta_{i+1}(q^{s})}{\theta_{i+1}(q^{i+1})}
    =
    (q^{r+s-i}-q^{r})
    \frac{\theta_{i}(q^{r})\theta_{i}(q^{s})}{\theta_{i}(q^{i})}
    +
    \frac{\theta_{i+1}(q^{r})\theta_{i+1}(q^{s})}{\theta_{i+1}(q^{i+1})}.
    $$
    Using the identities
    $$
    \theta_{i+1}(q^{s})=\theta_{i}(q^{s})(q^{s}-q^{i})
    \quad\text{and}\quad
    \theta_{i+1}(q^{r+1})=\theta_{i}(q^{r})(q^{r+i+1}-q^{i}),
    $$
    this is easy to check.  
\end{proof}

This proposition clarifies how it is that infinite sums can be
multiplied without producing infinite coefficients:
$$
\left(\sum_{r\ge0}a_{r}\phi_{r}\right)
\left(\sum_{s\ge 0}b_{s}\phi_{s}\right)
=
\sum_{j\ge 0}
\left(
\sum_{\substack{r,s\le j\\ r+s\ge j}} 
a_{r}b_{s}
\frac{\theta_{r+s-j}(q^{r})\theta_{r+s-j}(q^{s})}{\theta_{r+s-j}(q^{r+s-j})}
\right)
\phi_{j},
$$
where the important point is that the inner summations are finite.

The multiplicative structure of $\kuzk_{(p)}$ is very intricate; we
study it in more detail in Section~\ref{ringstructure}.

The following generalisation of Proposition~\ref{phi-prod} is used in
Section~\ref{ringstructure}.

\begin{prop}\label{phiprodexp}
    For all $r,s\ge m$,
    $$
    \phi_{r}\phi_{s}
    =
    \phi_{m}
    \!\!\sum_{j=\max(r,s)}^{r+s-m} c_{r,s}^{m,j}\phi_{j},
    $$
    where the coefficients are given by
    $$
    c_{r,s}^{m,j} 
    =
    \frac{\theta_{r+s-j}(q^{r})\theta_{r+s-j}(q^{s})}
    {q^{(r+s-j-m)m}\theta_{r+s-j-m}(q^{r+s-j-m})
    \theta_{m}(q^{r})\theta_{m}(q^{s})}.
    $$
\end{prop}

\begin{proof}
    Interchanging the role of~$r$ and~$s$ if necessary, we may assume
    that $s\ge r$.  Then Proposition~\ref{prodexp} provides the given
    expansion for $\phi_{r}\phi_{s}$ with
    $c_{r,s}^{m,j}=A_{r-m,j-s}(\bq[m],\bq[s])$.  The formula for
    $c_{r,s}^{m,j}$ holds since this expression satisfies the
    recurrence~\eqref{prodrecur}.
\end{proof}

The multiplication formula for the $f_{n}(w)$ (Proposition~7
of~\cite{ccw}) leads by duality to the following result.

\begin{prop}\label{coproduct}
    The coproduct satisfies
    $$
    \Delta\phi_{n}
    =
    \sum_{\substack{r,s\ge 0\\ r+s\le n}}
    \frac
    {\theta_{r+s}(q^{n})}
    {\theta_{r}(q^{r})\theta_{s}(q^{s})}
    \;
    \phi_{n-r}\tensor\phi_{n-s}.
    $$
    \qed
\end{prop}

The bounds in this summation ensure that the formula determines a map
from $\kuzk_{(p)}$ into the completed tensor product
$\kuzk_{(p)}\widehat\tensor\kuzk_{(p)}$, whose elements may be written
as doubly infinite sums of the~$\phi_{i}\tensor\phi_{j}$.  It is
unreasonable to expect a coproduct to map into the usual tensor
product in this setting, so we will take the terms coalgebra,
bialgebra and Hopf algebra to mean structures where the coproduct maps
into the completed tensor product.

\section{The ring structure of $k^{0}(k)_{(p)}$}
\label{ringstructure}

In this section we write~$I$ for the augmentation ideal of the
algebra $\kuzk_{(p)}$ of stable operations of degree~0 in \plocal\
connective \Ktheory.  Thus $I/I^{2}$ is the module of indecomposables. 
The ring of \padic\ integers is denoted by~\Zp.

\begin{thm}\label{ringstructurethm}\quad\par
    \begin{enumerate}
        \item\label{rs1}
        The ring $\kuzk_{(p)}$ is not Noetherian.
    
        \item\label{rs2}
        Its module of indecomposables is isomorphic to~\Zp.
    \end{enumerate}
\end{thm}

To prove this theorem, we define a sequence of ideals of~$\kuzk_{(p)}$.

\begin{defn}
    For $m\ge 0$ let
    $$
    B_{m}
    =
    \bigset{\sum_{n\ge m}a_{n}\phi_{n}:a_{n}\in\Zpl}.
    $$
\end{defn}
     
Thus $B_{0}\supset B_{1}\supset B_{2}\supset\cdots$, with $B_{0}=\kuzk_{(p)}$
and $B_{1}=I$.  It is clear from Remark~\ref{coeff-action} that 
$B_{m}$ consists of those operations which act as zero on the 
coefficient groups~$\pi_{2i}(k_{(p)})$ for $0\le i<m$.  Thus this 
filtration is independent of our choice of~$q$.

\begin{prop}\label{Bprod}
    For each $0\le n\le m$, we have $B_{n}B_{m}=\phi_{n}B_{m}$.  In
    particular, $B_{m}$ is an ideal of~$\kuzk_{(p)}$.
\end{prop}

\begin{proof}
    Clearly $\phi_{n}B_{m}\subseteq B_{n}B_{m}$, since
    $\phi_{n}\in B_{n}$.

    Suppose that $\alpha=\sum_{r\ge n}a_{r}\phi_{r}\in B_{n}$ and
    $\beta=\sum_{s\ge m}b_{s}\phi_{s}\in B_{m}$, then, using
    Proposition~\ref{phiprodexp},
    \begin{align*}
        \alpha\beta
        &=
        \sum_{\substack{r\ge n\\s\ge m}}a_{r}b_{s}\phi_{r}\phi_{s}\\
        &=
        \sum_{\substack{r\ge n\\s\ge m}}a_{r}b_{s}
        \phi_{n}\!\!\sum_{j=\max(r,s)}^{r+s-n} c_{r,s}^{n,j}\phi_{j}\\
        &=
        \phi_{n}\sum_{j\ge m}
        \biggl(
        \sum_{\substack{r,s\le j\\r+s\ge j+n}}a_{r}b_{s} c_{r,s}^{n,j}
        \biggr)
        \phi_{j},
    \end{align*}
    which belongs to $\phi_{n}B_{m}$, as the inner summation is finite.
\end{proof}

\begin{lemma}\label{quotlemma}
    For each $m\ge 1$, the quotient $B_{m}/\phi_{1}B_{m}$ is
    isomorphic to~\Zp.
\end{lemma}

\begin{proof}
    Define a $\Zpl$-module homomorphism $\pi_{m}:B_{m}\to\Zp$ by
    $$
    \pi_{m}:
    \left(\sum_{n\ge m}a_{n}\phi_{n}\right)
    \mapsto
    \sum_{n\ge m}a_{n}(1-q^{m})(1-q^{m+1})\dots(1-q^{n-1}).
    $$
    Note that the sum does indeed converge \padically.

    By Proposition~\ref{phi-prod},
    $\phi_{1}\phi_{n}=(q^{n}-1)\phi_{n}+\phi_{n+1}$.  Hence
    $\pi_{m}(\phi_{1}\phi_{n})=0$, if $n\ge m$, and
    $\ker\pi_{m}\supseteq\phi_{1}B_{m}$.

    Now suppose $\alpha=\sum_{n\ge m}a_{n}\phi_{n}\in\ker\pi_{m}$.  We
    will show that there exists 
    $\beta=\sum_{n\ge m}b_{n}\phi_{n}\in B_{m}$ such that
    $\alpha=\phi_{1}\beta$.  This is equivalent to showing that the
    equations
    \begin{equation}\label{beqns}
        b_{n-1}+(q^{n}-1)b_{n}
        =
        a_{n}
        \qquad(n\ge m)
    \end{equation}
    may be solved for \set{b_{n}\in\Zpl:n\ge m}, where $b_{n}=0$ for
    $n<m$.

    Suppose that we have found $b_{r}\in\Zpl$ satisfying~\eqref{beqns}
    for $m\le r<n$.  It follows then that
    \begin{align*}
        0=\pi_{m}(\alpha)
        &=
        (a_{n}-b_{n-1})(1-q^{m})\dots(1-q^{n-1})\\
        &\qquad
        +a_{n+1}(1-q^{m})\dots(1-q^{n})+\dotsb.
    \end{align*}
    Thus if $M=\nu_{p}\bigl((1-q^{m})\dots(1-q^{n-1})\bigr)$ and
    $N=\nu_{p}(1-q^{n})$,
    $$
    0
    \equiv
    (a_{n}-b_{n-1})p^{M}\mod{p^{M+N}},
    $$
    so that $a_{n}\equiv b_{n-1}\mod{p^{N}}$, and \eqref{beqns} may be
    solved for $b_{n}\in\Zpl$.

    This shows that $\ker\pi_{m}=\phi_{1}B_{m}$; it remains to prove
    that $\pi_{m}$ is surjective.

    Let $M=\nu_{p}\bigl((1-q)\dots(1-q^{m-1})\bigr)$ and choose $R$
    such that $m\le p^{R-1}(p-1)$.  For any $r\ge R$ we can write
    $$
    \prod_{j=m}^{p^{r-1}(p-1)}(1-q^{j})
    =
    u_{r}p^{\frac{p^{r}-1}{p-1}-M},
    $$
    where $u_{r}\in\Zplu$.

    If $x\in\Zp$, by grouping the \padic\ digits of~$x$ we can write
    $$
    x
    =
    x_{0}
    +x_{R}p^{\frac{p^{R}-1}{p-1}-M}
    +x_{R+1}p^{\frac{p^{R+1}-1}{p-1}-M}
    +\dots
    +x_{r}p^{\frac{p^{r}-1}{p-1}-M}
    +\dotsb,
    $$
    where $x_{r}\in\N$ for $r=0$, and for $r\ge R$.  Now let
    $$
    \alpha
    =x_{0}\phi_{m}
    +\sum_{r\ge R}\frac{x_{r}}{u_{r}}\phi_{p^{r-1}(p-1)+1}.
    $$
    It is clear that $\alpha\in B_{m}$ and $\pi_{m}(\alpha)=x$.
\end{proof}

\begin{prop}
    For $m\ge 1$, the ideal $B_{m}$ is not finitely generated.
\end{prop}

\begin{proof}
    By Proposition~\ref{phiprodexp}, $\kuzk_{(p)}$ acts on the quotient
    $B_{m}/\phi_{1}B_{m}=\Zp$ via the augmentation
    $\epsilon:\kuzk_{(p)}\to\Zpl$ and the inclusion $\Zpl\subset\Zp$. 
    Thus if $B_{m}$ were a finitely generated $\kuzk_{(p)}$-ideal,
    then $\Zp$ would be a finitely generated
    $\Zpl$\nobreakdash-module.  But a finite subset of~\Zp\ can only
    generate a countable \Zpl\nobreakdash-submodule of~\Zp.
\end{proof}

\begin{proof}[Proof of Theorem~\ref{ringstructurethm}]
    Part~\eqref{rs1} follows immediately.  For part~\eqref{rs2} we note that 
    Proposition~\ref{Bprod} shows that $I^{2}=\phi_{1}B_{1}$, so that 
    $I/I^{2}=B_{1}/\phi_{1}B_{1}\cong\Zp$ by Lemma~\ref{quotlemma}.
\end{proof}

It is interesting to note how far the augmentation ideal~$I$ is from
being generated by~$\phi_{1}=\Psi^{q}-1$.  Before carrying out the
completion giving $\kuzk_{(p)}$, the augmentation ideal
in~$\Zpl[\Psi^{q}-1]$ is principal, and it is so again after
$p$\nobreakdash-completion; see Section~\ref{pcomp}.  In contrast, we
have

\begin{cor}\label{Bbquot}
    The quotient $I/\langle \phi_{1} \rangle$ is isomorphic
    to~$\Zp/\Zpl$, where $\langle\phi_{1}\rangle$ is the ideal
    of~$\kuzk_{(p)}$ generated by~$\phi_{1}$.
\end{cor}

\begin{proof}
      It is clear that $\kuzk_{(p)}=\Zpl+I$.
\end{proof}

We remark that the abelian group $\Zp/\Zpl$ is torsion free and
divisible.  It is thus a $\Q$\nobreakdash-vector space.

\medskip

The intersection of the ideal~$B_{m}$ with the polynomial subalgebra
of~$\kuzk_{(p)}$ generated by~$\Psi^{q}$ is the principal ideal
generated by~$\phi_{m}=\theta_{m}(\Psi^{q})$.  Theorem~\ref{phibasis}
exhibits $\kuzk_{(p)}$ as the completion of $\Z_{(p)}[\Psi^{q}]$ with
respect to the filtration by these ideals.  But note that this
filtration is not multiplicative, in the sense
of~\cite{Northcott1968}, and hence there is no associated graded ring. 
It is for this reason that the completion fails to be Noetherian, and
contains zero divisors, as we see in Section~\ref{adamscon}.

Of course the Adams operations~$\Psi^{r}$, where $r$ is a \plocal\
unit, generate inside~$\kuzk_{(p)}$ a copy of the group
ring~$\Zpl[\Zplu]$ which is dense in the
$B_{m}$\nobreakdash-filtration topology.  Thus $\kuzk_{(p)}$ is
naturally a completion of~$\Zpl[\Zplu]$.

\medskip

In the following theorem we identify the units of~$\kuzk_{(p)}$ in 
terms of our basis.  This formulation is closely related to Theorem~1 
of~\cite{Johnson1987}.

\begin{thm}\label{kkunits}
    The element $\sum_{n\ge 0}a_{n}\phi_{n}$ is a unit in the ring
    $\kuzk_{(p)}$ if and only if
    $\sum_{n=0}^{i}a_{n}\theta_{n}(q^{i})$ is a \plocal\ unit for
    $i=0,1,\dots,p-2$.
\end{thm}

\begin{proof}
    If $\alpha=\sum_{n\ge 0}a_{n}\phi_{n}$ is a unit, then, by 
    Remark~\ref{coeff-action},
    $\sum_{n=0}^{i}a_{n}\theta_{n}(q^{i})$ represents the action 
    of~$\alpha$ on~$\pi_{2i}(k_{(p)})$, and so must be invertible for 
    all~$i$.
    
    Conversely, assume $\sum_{n=0}^{i}a_{n}\theta_{n}(q^{i})\in\Zplu$
    for $i=0,1,\dots,p-2$, then, since
    $\theta_{n}(q^{i})\equiv\theta_{n}(q^{j})\mod{p}$ if 
    $i\equiv j\mod{p-1}$ and $\theta_{n}(q^{i})=0$ if $n>i$, this holds
    for all~$i\ge 0$.

    Now suppose, inductively, that we have found
    $b_{0},b_{1},\dots,b_{i-1}\in\Z_{(p)}$ such that
    $$
    \left(\sum_{n\ge 0}a_{n}\phi_{n}\right)
    (b_{0}+b_{1}\phi_{1}+\dots+b_{i-1}\phi_{i-1})
    \in
    1+B_{i}.
    $$
    Then, using Proposition~\ref{phi-prod}, we see 
    that, for any $b_{i}$ in~\Zpl,
    $$
    \left(\sum_{n\ge 0}a_{n}\phi_{n}\right)
    (b_{0}+b_{1}\phi_{1}+\dots+b_{i}\phi_{i})
    \in
    1+B_{i},
    $$
    with the coefficient of~$\phi_{i}$ having the form
    $$
    b_{i}
    \left(\textstyle\sum_{n=0}^{i}a_{n}\theta_{n}(q^{i})\right)
    +\text{terms involving $b_{1},\dots,b_{i-1}$}.
    $$
    Thus, by our assumption, it is possible to chose $b_{i}\in\Zpl$ so that
    $$
    \left(\sum_{n\ge 0}a_{n}\phi_{n}\right)
    (b_{0}+b_{1}\phi_{1}+\dots+b_{i}\phi_{i})
    \in
    1+B_{i+1}.
    $$
    Repeating this process ad~infinitum yields the required inverse.
\end{proof}

\section{The Adams splitting and the Adams summand}
\label{adamscon}

Idempotent operations which split \plocal\ \Ktheory\ (for $p$~odd)
into $p-1$ summands were constructed in~\cite{Adams1969}.  We show how to
write the connective version of these idempotents in terms of our
basis elements.

\begin{prop}\label{conn-idem}
    If $\alpha\in\set{0,1,\dots,p-2}$, the Adams idempotent
    $e_{\alpha}\in\kuzk_{(p)}$ is given by
    $$
    e_{\alpha}
    =
    \sum_{n\ge 0}c_{n,\alpha}\phi_{n},
    $$
    where
    $$
    c_{n,\alpha} 
    = 
    \frac{1}{\theta_{n}(q^{n})}
    \sum(-1)^{n-i}q^{\binom{n-i}{2}}\gp{n}{i},
    $$
    the summation being over all integers~$i$ for which $0\le i\le n$
    and $i\equiv\alpha\mod{p-1}$.
\end{prop}

\begin{proof}
    If an operation~$\phi\in\kuzk_{(p)}$ acts on the coefficient
    group~$\pi_{2i}(k_{(p)})$ by multiplication by $\lambda_{i}$, then
    $\langle\phi,w^{i}\rangle=\lambda_{i}$.  Hence Proposition~8
    of~\cite{ccw} shows that
    $$
    \bigl\langle\phi,f_{n}(w)\bigr\rangle
    =
    \frac{1}{\theta_{n}(q^{n})}
    \sum_{i=0}^{n}(-1)^{n-i}q^{\binom{n-i}{2}}\gp{n}{i}\lambda_{i}.
    $$
    The result now follows by duality, since $e_{\alpha}$ acts as the
    identity on $\pi_{2i}(k_{(p)})$ if $i\equiv\alpha\mod{p-1}$, and
    as zero otherwise.
\end{proof}

The spectra $K_{(p)}$ and~$k_{(p)}$ are each split by Adams's
idempotents into $p-1$ suspensions of a multiplicative spectrum which
we denote by~$G$ and~$g$, respectively.%
\footnote{The notations $E(1)$~and~$e(1)$, or $L$~and~$l$, are also used.}
The coefficient ring $G_{*}$ can be identified with the subring
$\Zpl[u^{p-1},u^{-p+1}]\subset\Zpl[u,u^{-1}]=\pi_{*}(K_{(p)})$, and
$g_{*}$ is identified with
$\Zpl[u^{p-1}]\subset\Zpl[u]=\pi_{*}(k_{(p)})$.

There is an algebra isomorphism $\iota:\guzg\to e_{0}\kuzk_{(p)}$,
under which $\Psi^{q}\in\guzg$ maps to~$e_{0}\Psi^{q}$, but
$e_{0}\kuzk_{(p)}$ is not a sub-bialgebra of $\kuzk_{(p)}$.  However,
composing the projection $\kuzk_{(p)}\to e_{0}\kuzk_{(p)}$ with the
inverse of $\iota$, exhibits $\guzg$ as a quotient bialgebra
of~$\kuzk_{(p)}$.  Thus \guzg\ is a summand of~$\kuzk_{(p)}$ as an
algebra, but not as a bialgebra.

Note that if $p>3$, then, in contrast to the \padic\
case~\cite{Mitchell1993, Clarke1987}, the algebra $\kuzk_{(p)}$ is not
isomorphic to $\guzg\widehat{\otimes}\Z_{(p)}[C_{p-1}]$.  This can be
seen by considering the action on the coefficient ring which shows
that $\kuzk_{(p)}$ contains no elements of order~$p-1$.

\medskip

The results of Section~\ref{hopfalgcon} have analogues for the
algebra~$\guzg$ of degree zero stable operations on the Adams summand. 
We need first to adapt the ideas of~\cite{ccw} to provide a basis
for~\Gzg.  Let $z=w^{p-1}\in\Gzg$.

Recalling that we have chosen $q$ to be primitive modulo~$p^{2}$, we
write $\qh=q^{p-1}$ and let
$$
\thetah_{n}(X)
=
\prod_{i=0}^{n-1} (X-\qh^{i}).
$$

\begin{prop}
    A $\Zpl$-basis for \Gzg\ is given by the elements
    $$
    \fh_{n}(z)
    =
    \frac{\thetah_{n}(z)}{\thetah_{n}(\qh^{n})},
    \qquad\text{for $n\ge 0$.}
     $$
\end{prop}

\begin{proof}
    We have
    $$
    \Gzg 
    =
    \bigset{f(z)\in\Q[z]:f(1+p\Zpl)\subseteq\Zpl}.
    $$
    The multiplicative group $1+p\Zp$ is topologically generated
    by~$\qh$, and $\nu_{p}(\qh^{n}-1)=1+\nu_{p}(n)$ for all $n\ge 1$. 
    The rest of the proof parallels that of~\cite[Proposition~3]{ccw}.
\end{proof}

We now identify the dual basis for the algebra~$\guzg=\Guzg$.

\begin{defn}
    Define $\phih_{n}\in\guzg$, for $n\ge 0$, by
    $$
    \phih_{n}
    =
    \thetah_{n}(\Psi^{q})
    =
    \prod_{i=0}^{n-1}(\Psi^{q}-\qh^{i}).
    $$
    To avoid any misunderstanding, we emphasise that the~$q$ indexing
    the Adams operation is hatless.
\end{defn}

\begin{thm}[{\cite[Theorem~2.2]{Lellman1984}}]\label{hphibasis}
    The elements of \guzg\ can be expressed uniquely as infinite sums
    $$
    \sum_{n\ge 0}a_{n}\phih_{n},
    $$
    where $a_{n}\in\Z_{(p)}$.
\end{thm}

\begin{proof}
    The proof is just as for Theorem~\ref{phibasis}: $\guzg=G^{0}(g)$
    is $\Zpl$\nobreakdash-dual to the bialgebra~\Gzg, and the
    $\phih_{n}$ are dual to the $\fh_{n}(z)$, the action being given
    by $\phih_{m}\cdot\fh_{n}(z)=\qh^{m(m-n)}z^{m}\fh_{n-m}(z)$.
\end{proof}

The formulas for the product and coproduct of the~$\phih_{n}$ in~\guzg\
are, of course, just like those of Section~\ref{hopfalgcon} for
the~$\phi_{n}$, but with the primitive element~$q$ replaced by~$\qh$. 
Similarly the proof of Theorem~\ref{ringstructurethm} generalises to
show that \guzg\ is not Noetherian.

The proof of Theorem~\ref{kkunits} simplifies in the split context, since
$\thetah_{n}(\qh^{i})$ is divisible by~$p$ for all $i\ge 0$ and all
$n>0$.  Hence we have

\begin{thm}\label{ggunits}
    An element of~\guzg\ is a unit if and only if its augmentation is 
    a unit.
    \qed
\end{thm}

This result was proved by Johnson in~\cite{Johnson1987}; it shows that
\guzg\ is a local ring.  Johnson also showed in that paper that \guzg\
is an integral domain.

\section{Operations in $p$-adic $K$-theory}
\label{pcomp}

The \padic\ completion of \guzg, which we denote by $\guzg_{p}$, is 
the ring of degree zero operations in \padic\ \Ktheory.  In fact in 
this case the algebra of operations in the connective theory does not 
differ from the algebra of operations in the periodic theory; 
see~\cite{Johnson1986} and Section~\ref{convper} below.  We give here 
an algebraic proof of the result due to Clarke~\cite{Clarke1987} and 
Mitchell~\cite{Mitchell1993} that $\guzg_{p}$ is a power series 
ring.  

\begin{thm}
    $\guzg_{p}$ is equal to the power series ring over~\Zp\ generated
    by~$\Psi^{q}-1$.
\end{thm}

\begin{proof}
    Retaining the notation of Section~\ref{adamscon}, we let
    $s(n,i),S(n,i)\in\Zpl$ be such that, for $n\ge 1$,
    $$
    \thetah_{n}(X)
    =
    \sum_{i=1}^{n}s(n,i)(X-1)^{i},
    \quad\text{and}\quad
    (X-1)^{n}
    =
    \sum_{i=1}^{n}S(n,i)\thetah_{i}(X).
    $$
    These constants are analogues of the Stirling numbers of the first
    and second kinds, respectively.  (In the notation of
    Section~\ref{polyid}, $s(n,i)=A_{n,i}(\bqh,\bone)$ and
    $S(n,i)=A_{n,i}(\bone,\bqh)$, where \bqh\ is the sequence
    $(\qh^{i-1})_{i\ge 1}$ and \bone\ is the constant sequence
    $(1)_{i\ge 1}$.)  It is clear that $s(n,n)=S(n,n)=1$,
    $s(n,1)=(1-\qh)(1-\qh^{2})\dots(1-\qh^{n-1})$, and
    $S(n,1)=(\qh-1)^{n-1}$.  Moreover the two lower triangular
    matrices $\bigl(s(n,i)\bigr)_{n,i\ge 1}$ and
    $\bigl(S(n,i)\bigr)_{n,i\ge 1}$ are mutually inverse.

    In these cases, the recurrence of Proposition~\ref{basicrecur}
    becomes
    \begin{align*}
        s(n+1,i)
        &=s(n,i-1)-(\qh^{n}-1)s(n,i),\\
        \text{and}\quad
        S(n+1,i)
        &=S(n,i-1)+(\qh^{i}-1)S(n,i).
    \end{align*}
    Since $\qh^{n}-1$ is divisible by~$p$ for all~$n$, it follows easily
    that $\nu_{p}\bigl(s(n,i)\bigr)\ge n-i$ and
    $\nu_{p}\bigl(S(n,i)\bigr)\ge n-i$.

    Let $(p,Y)$ denote the maximal ideal of~$\Z_{p}[[Y]]$.  Since
    $\thetah_{n}(Y+1)=Y(Y+1-\qh)\dots(Y+1-\qh^{n-1})$, we have
    $\thetah_{n}(Y+1)\in(p,Y)^{n}$ for all~$n$.  As a result, there
    are ring homomorphisms forming the following commutative diagram
    $$
    \begin{CD}
        \Zp[X]/\bigl(\thetah_{n}(X)\bigr)
        @>>> \Zp[[Y]]/(p,Y)^{n}\\
        @AAA
        @AAA\\
        \Zp[X]/\bigl(\thetah_{n+1}(X)\bigr)
        @>>> \Zp[[Y]]/(p,Y)^{n+1}
    \end{CD}
    $$
    in which the horizontal maps are defined by $X\mapsto Y+1$.  In
    the limit there is a ring homomorphism
    $$
    \guzg_{p}
    =
    \invlim\Zp[X]/\bigl(\thetah_{n}(X)\bigr)
    \to
    \invlim\Zp[[Y]]/(p,Y)^{n}
    =
    \Zp[[Y]],
    $$
    which, we will show, is an isomorphism.  The variable~$X$
    corresponds to~$\Psi^{q}$, and thus $Y$ to~$\Psi^{q}-1$.

    The kernel of
    $\Zp[X]/\bigl(\thetah_{n}(X)\bigr)\to\Zp[[Y]]/(p,Y)^{n}$, which we
    shall denote by~$I_{n}$, is the free \Zp-module generated by the
    elements \bigset{\bigl[p^{n-i}(X-1)^{i}\bigr]:0\le i\le n-1},
    where $\bigl[g(X)\bigr]$ denotes the coset of~$g(X)$ in the
    quotient ring $\Zp[X]/\bigl(\thetah_{n}(X)\bigr)$.  Under the
    homomorphism $I_{n+1}\to I_{n}$, the element
    $\bigl[p^{n+1-i}(X-1)^{i}\bigr]$ maps to
    $p\bigl[p^{n-i}(X-1)^{i}\bigr]$, if $i<n$, but
    $\bigl[p(X-1)^{n}\bigr]$ maps to
    $$
    {} - \sum_{j=1}^{n-1}p\frac{s(n,j)}{p^{n-j}}\bigl[p^{n-j}(X-1)^{j}\bigr].
    $$
    Note here that $s(n,j)/p^{n-j}\in\Zpl$ by the remarks above.  This
    shows that the image of~$I_{n+1}$ lies in~$pI_{n}$, and hence,
    since no non-zero element of~$I_{n}$ is infinitely divisible
    by~$p$, that $\invlim I_{n}=0$.  Thus $\guzg_{p}$ maps injectively
    into~$\Zp[[Y]]$.

    To prove that it does so surjectively we define maps
    $\Zp[[Y]]\to\Zp[X]/\bigl(\thetah_{n}(X)\bigr)$ by
    $$
    \sum_{r\ge 0}c_{r}Y^{r}
    \mapsto c_{0}
    +
    \sum_{i=1}^{n-1}
    \left(\sum_{j=i}^{\infty}S(j,i)c_{j}\right)
    \bigl[\thetah_{i}(X)].
    $$
    It is here, of course, that we need to be working over~\Zp.  The
    convergence of the infinite series is guaranteed since $p^{j-i}$
    divides~$S(j,i)$.

    It will turn out that we have a ring homomorphism.  But at this
    point we need only to know that it is a homomorphism of
    \Zp-modules, and this is trivial.  It is also trivial that the
    maps factor through the projection
    $$
    \Zp[X]/\bigl(\thetah_{n+1}(X)\bigr)
    \to
    \Zp[X]/\bigl(\thetah_{n}(X)\bigr)
    $$
    and so define a \Zp-module homomorphism
    $$
    \Zp[[Y]]
    \to
    \invlim\Zp[X]/\bigl(\thetah_{n}(X)\bigr)
    =
    \guzg_{p}.
    $$
    To verify that this is the inverse of the ring homomorphism
    constructed earlier we need to verify that the composition
    $$
    \Zp[[Y]]
    \to
    \Zp[X]/\bigl(\thetah_{n}(X)\bigr)
    \to
    \Zp[[Y]]/(p,Y)^{n}
    $$
    coincides with the natural map.  But since 
    $S(j,i)\equiv 0\mod{p^{n}}$ for $j\ge i+n$, the composition sends
    $\sum_{r\ge 0}c_{r}Y^{r}$ to
    \begin{multline*}
        c_{0}
        +
        \sum_{i=1}^{n-1}
        \sum_{j=i}^{i+n-1}
        \sum_{k=1}^{i}S(j,i)s(i,k)c_{j}[Y^{k}]\\
        =
        c_{0}
        +
        \sum_{k=1}^{n-1}
        \sum_{j=k}^{2n-2}
        \sum_{i=\max(k,j-n+1)}^{\min(j,n-1)}S(j,i)s(i,k)c_{j}[Y^{k}],
    \end{multline*}
    where $[Y^{k}]$ denotes the coset of $Y^{k}$ in
    $\Zp[[Y]]/(p,Y)^{n}$, and we note that $[Y^{k}]=0$ for $k\ge n$.

    Now $[Y^{k}]$ has order~$p^{n-k}$ in~$\Zp[[Y]]/(p,Y)^{n}$, and
    $S(j,i)s(i,k)$ is divisible by~$p^{j-k}$.  This means that if
    $j\ge n$, then $S(j,i)s(i,k)[Y^{k}]=0$ and the image of
    $\sum_{r\ge 0}c_{r}Y^{r}$ is
      $$
      c_{0}
      +\sum_{k=1}^{n-1}
      \sum_{j=k}^{n-1}
      \sum_{i=k}^{j}S(j,i)s(i,k)c_{j}[Y^{k}]
      =c_{0}+
      \sum_{k=1}^{n-1}c_{k}[Y^{k}],
      $$
      since $\sum_{i=k}^{j}S(j,i)s(i,k)=\delta_{j,k}$.  This completes
      the proof.
\end{proof}

\section{Operations in periodic $K$-theory}
\label{hopfalgper}

We now turn our attention to the periodic case.  We let $\KuzK_{(p)}$
denote the algebra of degree zero operations in \plocal\ periodic
\Ktheory.  We show in this section how the results of
Sections~\ref{hopfalgcon}--\ref{adamscon} extend to this context. 
Here the degree zero operations determine all stable operations, and
$\KuzK_{(p)}$ is a Hopf algebra so we determine the antipode as well
as the other parts of the structure.

For each non-negative integer~$n$, we define the polynomial
$\Theta_{n}(x)$ by
$$
\Theta_{n}(X)
=
\prod_{i=1}^{n}(X-\qb_{i}),
$$
where $\qb_{i}$ is the $i$\nobreakdash-th term of the sequence
$$
\bqb
=
(1,q,q^{-1},q^{2},q^{-2},q^{3},q^{-3},q^{4},\dotsc)
=
\bigl(q^{(-1)^{i}\floor{i/2}}\bigr)_{i\ge 1},
$$
i.e., $\Theta_{n}(X)=\theta_{n}(X;\bqb)$ in the notation of the 
appendix.

\begin{defn}
    Define elements $\Phi_{n}\in\KuzK_{(p)}$, for $n \ge 0$, by
    $$
    \Phi_{n}
    =
    \Theta_{n}(\Psi^{q}).
    $$
\end{defn}

Thus, for example, $\Phi_{0}=1$, $\Phi_{1}=\Psi^{q}-1$,
$\Phi_{2}=(\Psi^{q}-1)(\Psi^{q}-q)$ and
$\Phi_{3}=(\Psi^{q}-1)(\Psi^{q}-q)(\Psi^{q}-q^{-1})$.

\begin{thm}\label{Phibasis}
    The elements of $\KuzK_{(p)}$ can be expressed uniquely as infinite sums
    $$
    \sum_{n \ge 0}a_{n}\Phi_{n},
    $$
    where $a_{n}\in\Z_{(p)}$.
\end{thm}

\begin{proof}
    The proof is analogous to that of Theorem~\ref{phibasis}.  The
    polynomials $F_{n}(w) = w^{-\floor{n/2}}f_{n}(w)$ form a
    $\Zpl$\nobreakdash-basis for $\KzK_{(p)}$ according to Corollary~6
    of~\cite{ccw}, and the $\Phi_{n}$ are, modulo multiplication by
    units, dual to this basis.  In fact the Kronecker pairing
    satisfies
    $$
    \bigl\langle\Phi_{n},F_{j}(w)\bigr\rangle
    = 
    \begin{cases}
        q^{-n\floor{n/2}},& \text{if $n=j$,}\\
        0,& \text{otherwise.} 
    \end{cases}
    $$
    As in the proof of Theorem~\ref{phibasis}, this follows from a
    study of the action of the dual on $\KzK_{(p)}$, but the details
    are a little more complicated and are given in the following
    result.
\end{proof}

\begin{lemma}\label{periodicaction}\quad\par
    \begin{enumerate}
        \item\label{j<n}
        $\Phi_{n}\cdot F_{j}(w)=0$ if $j<n$;
         
        \item\label{j=n}
        $\displaystyle
        \Phi_{n}\cdot F_{n}(w)
        =
        \begin{cases}
            q^{-nk}w^{-k}, & \text{if $n=2k$,}\\
            q^{-nk}w^{k+1}, & \text{if $n=2k+1$;}
        \end{cases}
        $
        \item\label{j>n}
        $\Phi_{n}\cdot F_{j}(w)$ is divisible by $f_{j-n}(w)$ for
        $j>n$.
    \end{enumerate}
\end{lemma}

\begin{proof}
    \eqref{j<n}\quad 
    Recall that $\Psi^{q}$ acts as $\Psi^{q}\cdot f(w)=f(qw)$.  If
    $j<n$ and $-\floor{j/2}\le i\le\ceiling{j/2}$, then
    $\Psi^{q}-q^{i}$ is a factor of~$\Phi_{n}$, so that $\Phi_{n}\cdot
    w^{i}=0$.  But $F_{j}(w)$ is a Laurent polynomial of
    codegree~$-\floor{j/2}$ and degree $j-\floor{j/2}=\ceiling{j/2}$.

     \eqref{j=n}\quad 
     By the proof of~\eqref{j<n}, all monomials occurring in
     $F_{n}(w)$ are annihilated by~$\Phi_{n}$ except one.  If $n=2k$,
     this is the lowest degree monomial~$w^{-k}$, whose coefficient is
     $f_{2k}(0)=q^{\binom{2k}{2}}/\theta_{2k}(q^{2k})$, and we have
     \begin{align*}
         \Phi_{2k}\cdot F_{2k}(w)
         &=
         \Phi_{2k}\cdot\bigl(f_{2k}(0)w^{-k}\bigr)\\
         &=
         q^{\binom{2k}{2}}\frac{\Theta_{2k}(q^{-k})}{\theta_{2k}(q^{2k})}w^{-k}\\
         &=
         q^{-2k^{2}}w^{-k}.
     \end{align*}
     If $n=2k+1$, it is the highest degree monomial~$w^{k+1}$ which must
     be considered.  The leading coefficient is
     $1/\theta_{2k+1}(q^{2k+1})$, and
     \begin{align*}
         \Phi_{2k+1}\cdot F_{2k+1}(w)
         &=
         \frac{\Theta_{2k+1}(q^{k+1})}{\theta_{2k+1}(q^{2k+1})}w^{k+1}\\
         &=
         q^{-(2k+1)k}w^{k+1}.
     \end{align*}

     \eqref{j>n}\quad 
     The proof is by (finite) induction on~$n$.  Note that
     $\Phi_{0}\cdot F_{j}(w)=F_{j}(w)$ is certainly divisible
     by~$f_{j}(w)$.  Now assume that 
     $\Phi_{n}\cdot F_{j}(w)=f_{j-n}(w)G_{n,j}(w)$ for some Laurent
     polynomial $G_{n,j}(w)$.  Then, since
     $\Phi_{n+1}=(\Psi^{q}-q^{i})\Phi_{n}$ for some~$i$,
     \begin{align*}
         \Phi_{n+1}\cdot F_{j}(w)
         &=
         (\Psi^{q}-q^{i})\cdot f_{j-n}(w)G_{n,j}(w)\\
         &=
         f_{j-n}(qw)G_{n,j}(qw)-q^{i}f_{j-n}(w)G_{n,j}(w).
     \end{align*}
     But this is zero for $w=1,q,q^{2},\dots,q^{j-n-2}$ and therefore
     divisible by~$f_{j-n-1}(w)$.
\end{proof}

\begin{remark}
    The action of the infinite sum $\sum_{n\ge 0}a_{n}\Phi_{n}$ on the
    coefficient group $\pi_{2i}(K_{(p})$ is multiplication by the finite
    sum
    $$
    \sum_{n=0}^{2|i|}a_{n}\Theta_{n}(q^{i}) .
    $$
    Clearly the augmentation sends $\sum_{n\ge 0}a_{n}\Phi_{n}$ to~$a_{0}$.
\end{remark}

We have the following analogue of Proposition~\ref{phi_n-exp}.

\begin{prop}\label{Phi_n-exp}
    For all $n\ge 0$,
    \begin{equation*}
        \Phi_{n}
        =
        \sum_{j=0}^{n}(-1)^{n-j}q^{e(n,j)}\gp{n}{j}\Psi^{q^{j}},
    \end{equation*}
    where
    \begin{equation*}
        e(n,j)
        =
        \begin{cases}
            -(n-j)(j-1)/2,&
            \text{if $n$ is even,}\\
            -(n-j)j/2,&
            \text{if $n$ is odd.}
        \end{cases}
    \end{equation*}
\end{prop}

\begin{proof}
    Assuming, inductively, that the result holds for even $n=2k$,
    since $\Phi_{2k+1}=\Phi_{2k}(\Psi^{q}-q^{-k})$, the coefficient
    of~$\Psi^{q^{j}}$ in~$\Phi_{2k+1}$ is
    \begin{multline*}
        (-1)^{j+1}
        \left(
        q^{-(2k-j+1)(j-2)/2}\gp{2k}{j-1}+q^{-k-(2k-j)(j-1)/2}\gp{2k}{j}
        \right)\\
        \begin{aligned}
            &=
            (-1)^{j+1}q^{-(2k+1-j)j/2}
            \left(
            q^{2k-j+1}\gp{2k}{j-1}+\gp{2k}{j}
            \right)\\
            &=
            (-1)^{j+1}q^{-(2k+1-j)j/2}\gp{2k+1}{j}.
            \end{aligned}
    \end{multline*}
    The argument to show that the odd case implies the next even case 
    is similar.
\end{proof}

Conversely, the proof of Theorem~\ref{Phibasis} yields

\begin{prop}\label{Psi_j-exp2}
    If $j\in\Zplu$,
    $$
    \Psi^{j}
    =
    \sum_{n\ge 0}
    q^{n\floor{n/2}}j^{-\floor{n/2}}
    \frac{\theta_{n}(j)}{\theta_{n}(q^{n})}\Phi_{n}.
    $$
    In particular, for $i\in\Z$,
    $$
    \Psi^{q^{i}}
    =
    \sum_{n\ge 0}q^{(n-i)\floor{n/2}}\gp{i}{n}\Phi_{n}.
    $$
    Note that this is a finite sum for $i\ge 0$.
    \qed
\end{prop}

We now consider the antipode~$\chi$ of the Hopf algebra~$\KuzK_{(p)}$. 
In the dual~$\KzK_{(p)}$ the antipode is given by $w\mapsto w^{-1}$
(see~\cite{ahs}), while $\chi\Psi^{j}=\Psi^{j^{-1}}$.

\begin{prop}\label{antipode}
    The antipode in~$\KuzK_{(p)}$ is determined by
    $$
    \chi\Phi_{n}
    =
    \sum_{j\ge 2\floor{(n-1)/2}+1}\!
    \left(
    \sum_{i=0}^{n}
    (-1)^{n-i}q^{(i+j)\floor{j/2}+e(n,i)}
    \gp{n}{i}\gp{-i}{j}
    \right)
    \Phi_{j},
    $$
    where $e(n,i)$ is defined in Proposition~\ref{Phi_n-exp}.
\end{prop}

\begin{proof}
    Proposition~\ref{Phi_n-exp} shows that
    $$
    \chi\Phi_{n}
    =
    \sum_{i=0}^{n}(-1)^{n-i}\gp{n}{i}q^{e(n,i)}\Psi^{q^{-i}},
    $$
    and, by Proposition~\ref{Psi_j-exp2},
    $$
    \Psi^{q^{-i}}
    =
    \sum_{j\ge 0}q^{(i+j)\floor{j/2}}\gp{-i}{j}\Phi_{j}.
    $$
    The required formula, with summation over $j\ge 0$, now follows by
    substitution.  
    
    To see that in fact the coefficients are zero until
    $j=2\floor{(n-1)/2}+1$ we note that by duality the expansion of
    $\chi\Phi_{n}$ can also be obtained from expressing
    $F_{j}(w^{-1})$ as a linear combination of the~$F_{n}(w)$.  The
    remarks in the proof of Lemma~\ref{periodicaction}~\eqref{j<n}
    show that the Laurent polynomial $F_{j}(w^{-1})$ can be written as
    a \Zpl-linear combination of the~$F_{n}(w)$ for $n\le j+1$ if $j$
    is odd, and for $n\le j$ if $j$ is even.  Thus we may
    write~$\chi\Phi_{n}$ is an infinite linear combination of the
    $\Phi_{j}$ for $j\ge n$ if $n$ is odd, and for $j\ge n-1$ if $n$
    is even.
\end{proof}

The same approach, using Propositions~\ref{Phi_n-exp}
and~\ref{Psi_j-exp2}, yields the following formulas for the product
and coproduct.

\begin{prop}\label{Phiprod}
     For all $r,s\ge 0$,
     $$
     \Phi_{r}\Phi_{s}
     =
     \sum_{k=\max(r,s)}^{r+s}A_{r,s}^{k}\Phi_{k},
     $$
     where
     $$
     A_{r,s}^{k}
     =
     \sum_{i=0}^{r}\sum_{j=0}^{s}
     (-1)^{r+s-i-j}q^{e(r,i)+e(s,j)+(k-i-j)\floor{k/2}}
     \gp{r}{i}\gp{s}{j}\gp{i+j}{k}.
     $$
     \qed
\end{prop}

\begin{prop}\label{Phi-coproduct}
    The coproduct in~$\KzK_{(p)}$ satisfies
    \begin{multline*}
        \Delta\Phi_{n}
        =\\
        \begin{cases}
            \displaystyle
            \sum_{\substack{r,s\ge 0\\ r+s\le n}}
            C_{n}^{r,s}\;\Phi_{n-r}\tensor\Phi_{n-s},&
            \text{if $n$ is even,}\\[3ex]
            \displaystyle
            \sum_{\substack{r,s\ge 0\\ r+s\le n}}
            C_{n}^{r,s}\;\Phi_{n-r}\tensor\Phi_{n-s}
            +
            \sum_{\substack{r,s\ge 2\\ r+s=n+1}}
            C_{n}^{r,s}\;\Phi_{n-r}\tensor\Phi_{n-s},&
            \text{if $n$ is odd,}
        \end{cases}
    \end{multline*}
    where
    \begin{multline*}
        C_{n}^{r,s}
        =\\
        \sum_{k=0}^{\min(r,s)}
        (-1)^{k}q^{e(n,n-k)+(k-r)\floor{(n-r)/2}+(k-s)\floor{(n-s)/2}}
        \gp{n}{k}\gp{n-k}{n-r}\gp{n-k}{n-s}.
    \end{multline*}
    \qed
\end{prop}

\medskip

The results of Section~\ref{ringstructure} also apply in the periodic 
case.

\begin{thm}\label{ringperthm}\quad\par
    \begin{enumerate}
        \item\label{per-rs1}
        The ring $\KuzK_{(p)}$ is not Noetherian.
    
        \item\label{per-rs2}
        Its module of indecomposables is isomorphic to~\Zp.
    \end{enumerate}
\end{thm}

\begin{proof}
    The proof exactly parallels that of
    Theorem~\ref{ringstructurethm}.  We define the family of
    ideals of~$\KuzK_{(p)}$
    $$
    A_{m}
    =
    \bigset{\sum_{n\ge m}a_{n}\Phi_{n}:a_{n}\in\Zpl}.
    $$
    Then we have $A_{n}A_{m}=\Phi_{n}A_{m}$ if $0\le n\le m$.  Just as is
    Section~\ref{ringstructure}, we show that, for $m\ge 1$,
    $A_{m}/\Phi_{1}A_{m}\cong\Zp$, and thus $A_{m}$ is not finitely
    generated.
\end{proof}

We can give criteria analogous to those of Theorem~\ref{kkunits} for 
when a general element of~$\KuzK_{(p)}$ is a unit, but we omit the 
details.

\medskip

We can now generalise the results of Section~\ref{adamscon} to the
periodic case.

\begin{prop}\label{per-idem}
    If $\alpha\in\set{0,1,\dots,p-2}$, the Adams idempotent
    $E_{\alpha}\in\KuzK_{(p)}$ is given by
    $$
    E_{\alpha}
    =
    \sum_{n\ge 0}C_{n,\alpha}\Phi_{n},
    $$
    where
    $$
    C_{n,\alpha} 
    = 
    \frac{1}{\theta_{n}(q^{n})}
    \sum(-1)^{\ceiling{n/2}-i}q^{n\floor{n/2}+\binom{\ceiling{n/2}-i}{2}}
    \gp{n}{\floor{n/2}+i},
    $$
    the summation being over all integers~$i$ for which $-n-1<2i\le n+1$
    and $i\equiv\alpha\mod{p-1}$.
    \qed
\end{prop}

Recall that we let $\qh=q^{p-1}$, where $q$ is primitive modulo~$p^{2}$.

\begin{defn}
    Let
    $$
    \hat\Phi_{n}
    = 
    \prod_{i=1}^{n}(\Psi^{q}-\qh^{(-1)^{i}\floor{i/2}}).
    $$
\end{defn}

\begin{thm}\label{hPhibasis}
    The elements of \GuzG\ can be expressed uniquely as infinite sums
    $$
    \sum_{n\ge 0}a_{n}\hat\Phi_{n},
    $$
    where $a_{n}\in\Z_{(p)}$.
    \qed
\end{thm}

It is possible to write down formulas generalising those above for the
antipode, product and coproduct in~\GuzG, but we omit the details.  
We can also show easily that \GuzG\ is a non-Noetherian local ring.

\section{The relation between the connective and periodic cases}
\label{convper}

We consider now the relation between the connective and periodic
cases, focussing on the non-split setting, although it is clear that
similar results hold for the relation between \guzg\ and~\GuzG.

The covering map $k \to K$ leads to an inclusion
$$
\KuzK_{(p)}
\hookrightarrow 
\Kuzk_{(p)}
=
\kuzk_{(p)},
$$
which is described by the following formula.

\begin{prop}\label{Phi-phi}
    For $n>0$,
    $$
    \Phi_{n} 
    = 
    \sum_{i=\floor{n/2}+1}^{n}
    \left(
    \sum_{j=i}^{n}(-1)^{n-j}q^{e(n,j)}\gp{n}{j}\gp{j}{i}
    \right)
    \phi_{i},
    $$
    where $e(n,j)$ is as defined in Proposition~\ref{Phi_n-exp}.
\end{prop}

\begin{proof}
    We obtain the required formula by combining 
    Propositions~\ref{Phi_n-exp} and~\ref{Psi_j-exp}, but with a 
    summation from~$i=0$ to~$n$.
    
    To see that the coefficients are zero for $i=0,\dots,\floor{n/2}$,
    we note that the polynomial~$\Theta_{n}(X)$ is divisible
    by~$\theta_{\floor{n/2}+1}(X)$, and the quotient is
    $\theta_{\ceiling{n/2}-1}\bigl(X,(q^{-1},q^{-2},\dotsc)\bigr)$, in
    the notation of~Section~\ref{polyid}.  Writing this quotient as a
    linear combination of the 
    $\theta_{r}\bigl(X,
    (q^{\floor{n/2}+1},q^{\floor{n/2}+2},\dotsc)\bigr)$, 
    and substituting $X=\Psi^{q}$, gives rise to the formula
   $$
   \Phi_{n} 
   = 
   \sum_{i=\floor{n/2}+1}^{n}
   A_{n,i-\floor{n/2}-1}
   \bigl(
   (q^{-1},q^{-2},\dotsc),
   (q^{\floor{n/2}+1},q^{\floor{n/2}+2},\dotsc)
   \bigr)
   \phi_{i}.
   $$
\end{proof}

An arbitrary element of $\KuzK_{(p)}$ can then be written as 
$$
\sum_{n\ge 0}a_{n}\Phi_{n}
=
\sum_{i\ge 0}
\left(
\sum_{n=i}^{2i}
\left(
\sum_{j=i}^{n}(-1)^{n-j}q^{e(n,j)}\gp{n}{j}\gp{j}{i}
\right)
a_{n}
\right)
\phi_{i}.
$$
Note that the inner summations are finite.

Since they are polynomials in~$\Psi^{q}$, the basis elements
$\phi_{i}$ lie in the image of the inclusion and can be expressed in
terms of our basis for~$\KuzK_{(p)}$ as follows.

\begin{prop}\label{phi-Phi}
    $$
    \phi_{n}
    =
    \sum_{i=0}^{n}
    q^{\floor{(i+3)/2}(n-i)}\theta_{n-i}(q^{-i-1})
    \gp{n-\ceiling{i/2}-1}{n-i}
    \Phi_{i}.
    $$
\end{prop}

\begin{proof}
    It is merely necessary to verify that the given coefficients 
    satisfy the recurrence given by Proposition~\ref{basicrecur} 
    for~$A_{n,i}(\bq,\bqb)$, and this is routine.
\end{proof}

Note that the summation here runs from $i=0$.  It is not the case that
the coefficients are zero for $i$ small enough.  So any attempt to
write an arbitrary infinite sum $\sum a_{n}\phi_{n}$ in~$\kuzk_{(p)}$
in terms of the~$\Phi_{i}$ will lead to infinite sums for the
coefficients.  This reflects the fact that the map
$\KuzK_{(p)}\hookrightarrow\kuzk_{(p)}$ is a {\em strict}
monomorphism.  However, if we complete \padically, the highly
$p$\nobreakdash-divisible factor $\theta_{n-i}(q^{-i-1})$ ensures that
these inner sums converge.  Thus we recover the fact that, by
contrast, $\KuzK_{p}\to\kuzk_{p}$ is an isomorphism, as discussed
in~\cite{Clarke1987} and~\cite{Johnson1986}.

\section{$2$-local operations}
\label{2-local}

Here we provide a description of operations in \twolocal\ \Ktheory. 
This is a little more complicated than for odd primes, essentially
because instead of the single `generator' $\Psi^{q}$ one has to deal
with both~$\Psi^{3}$ and~$\Psi^{-1}$.

Throughout this section and the next, the variable~$q$ occurring
implicitly in the polynomial~$\theta_{n}(X)$ will be set equal to~$9$,
and we let $\thetab_{n}(X)=\prod_{i=0}^{n-1}(X-3^{2i+1})$.  Thus, in
the notation of the appendix, $\theta_{n}(X)=\theta_{n}(X;\ba)$ and
$\thetab_{n}(X)=\theta_{n}(X;\bb)$, where \ba\ is the sequence of even
powers of~$3$ and \bb\ is the sequence of odd powers of~$3$.  These
choices are related to the fact that $\set{\pm 3^{i}:i\ge 0}$ is dense
in $\Ztu$; see~\cite{ccw}.

\begin{defn}
    Define elements $\zeta_{n}\in \kuzk_{(2)}$, for $n\ge 0$, by
    \begin{align*}
        \zeta_{2m+1}
        &=
        (\Psi^{-1}-1)\thetab_{m}(\Psi^{3}),\\
        \zeta_{2m}
        &=
        \theta_{m}(\Psi^{3})
        +
        \sum_{i=1}^{m}
        \frac{\theta_{i}(3)\theta_{i}(9^{m})}{2\theta_{i}(9^{i})}\zeta_{2m-2i+1}.
    \end{align*}
\end{defn}

\begin{thm}\label{2localops}
    The elements of $\kuzk_{(2)}$ can be expressed uniquely as
    infinite sums
    $$
    \sum_{n\ge 0}a_{n}\zeta_{n},
    $$
    where $a_{n}\in\Ztl$.
\end{thm}

\begin{proof}
    The proof mirrors that of Theorem~\ref{phibasis}.  We show that
    the given elements form the dual basis to the basis
    \set{f_{n}^{(2)}(w):n\ge 0} obtained for $\Kzk_{(2)}$
    in~\cite[Proposition~20]{ccw}.  This may be proved by induction
    arguments using the following formulas, describing the action
    of~$\kuzk_{(2)}$ on~$\Kzk_{(2)}$.
    \begin{align*}
        \Psi^{-1}\cdot f_{2m}^{(2)}(w)
        &=
        f_{2m}^{(2)}(w),\\
        \Psi^{-1}\cdot f_{2m+1}^{(2)}(w)
        &=
        f_{2m}^{(2)}(w)-f_{2m+1}^{(2)}(w),\\
        \Psi^{3}\cdot f_{2m}^{(2)}(w)
        &=
        9^{m}f_{2m}^{(2)}(w)+f_{2m-2}^{(2)}(w),\\
        \Psi^{3}\cdot f_{2m+1}^{(2)}(w)
        &=
        3^{2m+1}f_{2m+1}^{(2)}(w)-9^{m}f_{2m}^{(2)}(w)+f_{2m-1}^{(2)}(w).
    \end{align*}
\end{proof}

\begin{remark}
    The operation~$\zeta_{n}$ acts on the coefficient
    group~$\pi_{2i}(k_{(2)})$ as multiplication by the values given in
    the following table.
    
    \begin{table}[h]
        \centering
        \renewcommand{\arraystretch}{1.5}
        \begin{tabular}{c|cc}
            & $n=2m$& $n=2m+1$\\
            \hline
            $i$ even& $\theta_{m}(3^{i})$& $0$\\
            $i$ odd& $\thetab_{m}(3^{i})\strut$& $-2\thetab_{m}(3^{i})$
        \end{tabular}
    \end{table}
    
    \noindent
    Thus $\zeta_{n}$ acts as zero on $\pi_{2i}$ for all $i<n$.
\end{remark}

The following proposition gives product formulas.

\begin{prop}
    \begin{align*}
        \zeta_{2m}\zeta_{2n}
        &=
        \sum_{i=0}^{\min(m,n)}
        d_{m,n}^{i}
        \left(
        \zeta_{2m+2n-2i}-\tfrac{3^{i}-1}{2}\zeta_{2m+2n+1-2i}
        \right),\\
        \zeta_{2m+1}\zeta_{2n}
        &=
        \sum_{i=0}^{\min(m,n)}
        3^{i}d_{m,n}^{i}\zeta_{2m+2n+1-2i},\\
        \zeta_{2m+1}\zeta_{2n+1}
        &=
        -2\sum_{i=0}^{\min(m,n)}
        3^{i}d_{m,n}^{i}\zeta_{2m+2n+1-2i},
    \end{align*}
    where 
    $d_{m,n}^{i}=
    \dfrac{\theta_{i}(9^{m})\theta_{i}(9^{n})}{\theta_{i}(9^{i})}$.
\end{prop}

\begin{proof}
    These formulas are proved by long but straightforward induction
    arguments.
\end{proof}

Since $\zeta_{2n+1}=\zeta_{1}\zeta_{2n}$, the following proposition
completely determines the coproduct.

\begin{prop}
    \begin{align*}
        \Delta\zeta_{1}
        &=
        \zeta_{1}\tensor 1+\zeta_{1}\tensor\zeta_{1}+1\tensor\zeta_{1},\\
        \Delta\zeta_{2n}
        &=
        \sum_{\substack{r,s\ge 0\\ r+s\le n}}
        \frac
        {\theta_{r+s}(9^{n})}
        {\theta_{r}(9^{r})\theta_{s}(9^{s})}
        \;
        \zeta_{2n-2r}\tensor\zeta_{2n-2s}.
    \end{align*}
\end{prop}

\begin{proof}
    These formulas may be deduced from the product formula for the
    dual basis of cooperations just as in the proof of 
    Proposition~\ref{coproduct}.
\end{proof}

In principle, similar methods will give a description of the periodic
case~$\KuzK_{(2)}$.  However, this will be even more complicated, and
we omit the details.

\section{Operations in~$KO$-theory}
\label{KO}

We consider $ko$ and $KO$ localised at $p=2$.  (For $p$~odd, these
spectra split as $(p-1)/2$ copies of~$g$ and~$G$, respectively.) 
Proofs are omitted since the arguments are just the same as those
in~\cite{ccw} and in earlier sections of this paper.

Just as in Section~\ref{2-local}, the variable $q$ used implicitly in
the polynomials $\theta_{n}(X)$ and $\Theta_{n}(X)$ is set equal
to~$9$.

\begin{defn}\label{hdef}
    Let $x=w^{2}=u^{-2}v^{2}\in KO_{0}(ko)$; see~\cite{ahs}.  Let
    $$
    h_{n}(x)
    =
    \frac{\theta_{n}(x)}{\theta_{n}(9^{n})}.
    $$
\end{defn}

\begin{prop}\label{hbasis}
    \quad\par
    \begin{enumerate}
        \item \set{h_{n}(x):n\ge 0} is a \Ztl-basis for 
        $KO_{0}(ko)_{(2)}$. 
    
        \item \set{x^{-\floor{n/2}}h_{n}(x):n\ge 0} is a
        $\Ztl$\nobreakdash-basis for $KO_{0}(KO)_{(2)}$.
    \end{enumerate}
    \qed
\end{prop}

\begin{thm}\label{kobasis}
    \quad\par
    \begin{enumerate}
        \item The elements of $ko^{0}(ko)_{(2)}$ can be expressed
        uniquely as infinite sums
        $$
        \sum_{n \ge 0}a_{n}\theta_{n}(\Psi^{3}),
        $$
        where $a_{n}\in\Ztl$.
        \item The elements of $KO^{0}(KO)_{(2)}$ can be expressed
        uniquely as infinite sums
        $$
        \sum_{n \ge 0}a_{n}\Theta_{n}(\Psi^{3}),
        $$
        where $a_{n}\in\Ztl$.
    \end{enumerate}
    \qed
\end{thm}

\begin{cor}
    There is an isomorphism of bialgebras
    $$
    \kuzk_{(2)}/\langle\Psi^{-1}-1\rangle
    \cong
    ko^{0}(ko)_{(2)}.
    $$
\end{cor}

\begin{proof}
    We have seen that the ideal~$I$ generated by
    $\zeta_{1}=\Psi^{-1}-1$ is a coideal.  Note that
    $\zeta_{2n}\equiv\theta_{n}(\Psi^{3})$ modulo~$I$.  By
    Theorem~\ref{kobasis} the~$\theta_{n}(\Psi^{3})$ form a basis of
    $ko^{0}(ko)_{(2)}$, and the product and coproduct are respected.
\end{proof}

The methods of Section~\ref{pcomp} may be adapted to show that the
$2$\nobreakdash-adic completion $ko^0(ko)_2$ is the power series ring
over~$\Z_2$ generated by $\Psi^3-1$.

\section{Appendix: polynomial identities}
\label{polyid}

Given two sequences $\ba=(a_{i})_{i\ge 1}$ and $\bb=(b_{i})_{i\ge 1}$
of elements of a commutative ring~$R$, let
$$
\theta_{n}(X;\ba)
=
\prod_{i=1}^{n}(X-a_{i})
\quad\text{and}\quad
\theta_{n}(X;\bb)
=
\prod_{i=1}^{n}(X-b_{i}).
$$
These polynomials provide two bases for $R[X]$ as an
$R$\nobreakdash-module and can therefore be written in terms of each
other.

\begin{defn}\label{AnrDefn}
    Define $A_{n,r}(\ba,\bb)\in R$ by
    $$
    \theta_{n}(X;\ba)
    =
    \sum_{r=0}^{n}A_{n,r}(\ba,\bb)\theta_{r}(X;\bb).
    $$
\end{defn}

Together with the identities $A_{0,0}(\ba,\bb)=1$,
$A_{0,r}(\ba,\bb)=0$ for $r>0$, and $A_{n,-1}(\ba,\bb)=0$ for 
$n\ge 0$, the coefficients $A_{n,r}(\ba,\bb)$ are determined by the
following recurrence, which is used repeatedly in this paper.

\begin{prop}\label{basicrecur}
    For $n,r\ge 0$,
    \begin{equation*}
        A_{n+1,r}(\ba,\bb)
        =
        (b_{r+1}-a_{n+1})A_{n,r}(\ba,\bb)
        +
        A_{n,r-1}(\ba,\bb).
    \end{equation*}
\end{prop}

\begin{proof}
    Since $\theta_{n+1}(X;\ba)=(X-a_{n+1})\theta_{n}(X;\ba)$,
    \begin{align*}
        \theta_{n+1}(X;\ba)
        &=
        \sum_{r=0}^{n}A_{n,r}(\ba,\bb)(X-a_{n+1})\theta_{r}(X;\bb)\\
        &=
        \sum_{r=0}^{n}A_{n,r}(\ba,\bb)\bigl(
        \theta_{r+1}(X;\bb)+(b_{r+1}-a_{n+1})\theta_{r}(X;\bb)
        \bigr),
    \end{align*}
    and the result follows using
    $A_{n,-1}(\ba,\bb)=A_{n,n+1}(\ba,\bb)=0$.
\end{proof}

In most cases which we have considered, direct substitution in the
recurrence allows considerable simplification.  However it is possible
to give an explicit formula for the coefficients $A_{n,r}(\ba,\bb)$ as
polynomials in the~$a_{i}$ and~$b_{j}$.

\begin{prop}\label{AnrExplicit}
    For $n,r\ge 0$,
    $$
    A_{n,r}(\ba,\bb)
    =
    \sum_{\substack{J\subset\lbrace 1,\dots,n\rbrace\\|J|=r}}
    \prod_{\substack{1\le i\le n\\i\notin J}}
    (b_{\sigma(i,J)}-a_{i}),
    $$
    where $\sigma(i,J)=1+\bigl|\set{j\in J:j<i}\bigr|$.
\end{prop}

\begin{proof}
    Let
    $$
    \At_{n,r}
    =
    \sum_{\substack{J\subset\lbrace 1,\dots,n\rbrace\\|J|=r}}
    \prod_{\substack{1\le i\le n\\i\notin J}}
    (b_{\sigma(i,J)}-a_{i}).
    $$
    We have $\At_{0,0}=1$, $\At_{0,r}=0$ for $r>0$, and $\At_{n,-1}=0$
    for $n\ge 0$.  It is therefore only necessary to verify that
    $\At_{n,r}$ satisfies the recurrence
    $\At_{n+1,r}=(b_{r+1}-a_{n+1})\At_{n,r}+\At_{n,r-1}$ of
    Proposition~\ref{basicrecur}.
    
    Break the sum defining $\At_{n+1,r}$ into two parts by dividing
    the subsets $J\subset\lbrace 1,\dots,n+1\rbrace$ such that $|J|=r$
    according to whether $n+1\in J$ or not.  If $n+1\notin J$, then
    $J\subset\lbrace 1,\dots,n\rbrace$ and the corresponding summand
    in~$\At_{n,r}$ occurs in $\At_{n+1,r}$ multiplied by the factor
    $(b_{r+1}-a_{n+1})$ since $\sigma(n+1,J)=r+1$.

    If $n+1\in J$, let 
    $I=J\setminus\lbrace n+1\rbrace\subset\lbrace 1,\dots,n\rbrace$,
    then $\sigma(i,J)=\sigma(i,I)$ for all $i\notin J$, and
    $$
    \prod_{\substack{1\le i\le n+1\\i\notin J}}
    (b_{\sigma(i,J)}-a_{i})
    =
    \prod_{\substack{1\le i\le n\\i\notin I}}
    (b_{\sigma(i,I)}-a_{i})
    $$
    is a summand in both $\At_{n+1,r}$ and $\At_{n,r-1}$.

    It is clear that in both cases the process is reversible.
\end{proof}

We show finally how essentially the same coefficients arise in
formulas for products of the $\theta_{n}(X,\ba)$.  Given a sequence
$\bc=(c_{i})_{i\ge 1}$, we write~$\bc[m]$ for the shifted
sequence~$(c_{m+i})_{i\ge 1}$.

\begin{prop}\label{prodexp}
    If $r\ge m\ge 0$ and $s\ge 0$, then
    \begin{equation*}
        \theta_{r}(X;\bc)\theta_{s}(X;\bc)
        =
        \theta_{m}(X;\bc)
        \!\!
        \sum_{j=s}^{r+s-m}
        A_{r-m,j-s}(\bc[m],\bc[s])\theta_{j}(X;\bc).
    \end{equation*}
\end{prop}

\begin{proof}
    Setting $\ba=\bc[m]$ and $\bb=\bc[s]$ in Definition~\ref{AnrDefn},
    we have
    $$
    \theta_{r-m}(X;\bc[m])
    =
    \sum_{j=s}^{r+s-m}A_{r-m,j-s}(\bc[m],\bc[s])\theta_{j-s}(X;\bc[s]).
    $$
    Multiplying by $\theta_{m}(X;\bc)\theta_{s}(X;\bc)$, and using the
    identities
    $\theta_{r}(X;\bc)=\theta_{m}(X;\bc)\theta_{r-m}(X;\bc[m])$ and
    $\theta_{j}(X;\bc)=\theta_{s}(X;\bc)\theta_{j-s}(X;\bc[s])$, now
    gives the result.
\end{proof}

Writing $A_{r,s}^{m,j}=A_{r-m,j-s}(\bc[m],\bc[s])$ for the coefficient
in the above expansion for $\theta_{r}(X;\bc)\theta_{s}(X;\bc)$, the
recurrence of Proposition~\ref{basicrecur} takes the form
\begin{equation}\label{prodrecur}
    A_{r+1,s}^{m,j}
    =
    (c_{j+1}-c_{r+1})A_{r,s}^{m,j}+A_{r,s}^{m,j-1}.
\end{equation}

\end{document}